\documentclass[10pt,twoside]{article}
\usepackage{graphicx}
\usepackage{amscd, amssymb, amsthm, amsmath,eucal, verbatim}

\catcode`@=11
\def\ps@icm{\def\@oddhead{\hfill \leftmark \hfill\thepage }
\def\@evenhead{\thepage \hfill \rightmark \hfill}
\def\@oddfoot{}
\def\@evenfoot{}}
\def\ps@first{\def\@oddhead{ICM 2002 $\cdot$ Vol. ? $\cdot$ 1--10
\hfill}
\def\@evenhead{}
\def\@oddfoot{}
\def\@evenfoot{\hfill\copyright\ China Higher Education Press}}
\def\list#1#2{\ifnum \@listdepth >5\relax \@toodeep \else \global
\advance \@listdepth\@ne \fi \rightmargin \z@ \listparindent\z@
\itemindent\z@ \csname @list\romannumeral\the\@listdepth\endcsname
\def\@itemlabel{#1}\let\makelabel\@mklab \@nmbrlistfalse #2\relax
\@trivlist \parskip -\parsep \parindent\listparindent \advance
\linewidth -\rightmargin \advance\linewidth -\leftmargin \advance
\@totalleftmargin \leftmargin \parshape \@ne \@totalleftmargin
\linewidth \ignorespaces}
\catcode`@=12
\def\setzero{\setcounter{equation}{0}}
\def\thebibliography#1{\section*{References}
\list{[\arabic{enumi}]}{\settowidth \labelwidth{[#1]} \leftmargin
\labelwidth \advance \leftmargin \labelsep \usecounter{enumi}}
\def\newblock{\hskip .11em plus .33em minus .07em} \sloppy
\clubpenalty 4000 \widowpenalty 4000 \sfcode`\.=1000 \relax}
\newcommand{\proj}{\mathbf P}

\newcommand{\rarr}{\rightarrow}

\newcommand{\com}{\mathbb{C}}
\newcommand{\Q}{\mathbb{Q}}

\newcommand{\Z}{\mathbb{Z}}

\newcommand{\p}{\partial}

\def\scup{\mathbin{\text{\scriptsize$\cup$}}}

\newcommand{\LV}{\Lambda^{\frac\infty2}V}

\DeclareMathOperator{\ev}{ev}

\newcommand{\C}{\mathbb{C}}

\newtheorem{conj}{Conjecture}

\title{\bf  Three questions in Gromov-Witten theory\vskip 6mm}

\author{R. Pandharipande
\vspace*{-0.5cm}\thanks{Department of Mathematics,
Princeton University, Princeton, NJ 08544, USA. E-mail: rahulp@math.princeton.edu}}
\date{}

\begin{document}
\maketitle

\thispagestyle{first} \setcounter{page}{1}

\begin{abstract}\vskip 3mm
Three conjectural directions in Gromov-Witten theory are discussed:
Gorenstein properties, BPS states, and Virasoro constraints.
Each points to basic structures  in the subject which are not yet
understood.
\vskip 4.5mm


\end{abstract}

\vskip 12mm

\section{Introduction} \label{section 1}\setzero
\vskip-5mm \hspace{5mm }

Let $X$ be a nonsingular projective variety over $ \com$.
Gromov-Witten theory concerns integration over
$\overline{M}_{g,n}(X,\beta)$, the moduli space of stable maps
from genus $g$, $n$-pointed curves to $X$ representing the class
$\beta\in H_2(X,\Z)$. While substantial progress in the
mathematical study of Gromov-Witten theory has been made in the
past decade, several fundamental questions remain open. My goal
here is to describe three conjectural directions:
\begin{enumerate}
\item[(i)] Gorenstein properties of tautological rings,
\item[(ii)] BPS states for threefolds,
\item[(iii)] Virasoro constraints.
\end{enumerate}
Each points to basic structures in Gromov-Witten theory
which are not yet understood. New ideas in the subject will
be required for answers to these questions.

\section{Gorenstein properties of tautological rings} \label{tautor}
\setzero\vskip-5mm \hspace{5mm }

The study of the structure of the entire Chow ring of the
moduli space of pointed curves $\overline{M}_{g,n}$ appears quite
difficult at present.
As the principal motive is to understand
cycle classes obtained from algebro-geometric constructions,
we may restrict inquiry to the system of {\em tautological rings},
$R^*(\overline{M}_{g,n})$.
The tautological system is defined
to be the set of smallest $\Q$-subalgebras of the Chow rings,
$$R^*(\overline{M}_{g,n}) \subset A^*(\overline{M}_{g,n}),$$
satisfying the
following three properties:
\begin{enumerate}
\item[(i)] $R^*(\overline{M}_{g,n})$ contains
the cotangent line classes $\psi_1, \ldots, \psi_n$ where
$$\psi_i = c_1(L_i),$$
the first Chern class of the $i$th cotangent line bundle.
\item[(ii)] The system is closed under push-forward via
all maps forgetting markings:
$$\pi_*: R^*(\overline{M}_{g,n}) \rarr R^*(\overline{M}_{g,n-1}).$$
\item[(iii)] The system is closed under push-forward via
all gluing maps:
$$\pi_*: R^*(\overline{M}_{g_1,n_1\scup\{*\}})
\otimes_{\Q}
R^*(\overline{M}_{g_2,n_2\scup\{\bullet\}}) \rarr
R^*(\overline{M}_{g_1+g_2, n_1+n_2}),$$
$$\pi_*: R^*(\overline{M}_{g_1, n_1\scup\{*,\bullet\}}) \rarr
R^*(\overline{M}_{g_1+1, n_1}).$$
\end{enumerate}
Natural algebraic constructions typically yield Chow classes
lying in the tautological ring. See \cite{FP3}, \cite{GrP2} for further
discussion.

Consider the following basic
filtration of the moduli space of pointed curves:
$$\overline{M}_{g,n} \supset {M^c_{g,n}}
\supset M^{rt}_{g,n} \supset C_{g,n}.$$ Here, $M^c_{g,n}$ denotes
the moduli of pointed curves of compact type, $M^{rt}_{g,n}$
denotes the moduli of pointed curves with rational tails, and
$C_{g,n}$ denotes the moduli of pointed  curves with a fixed
stabilized complex structure $C_g$. The choice of $C_g$ will play
a role below.

The tautological rings for the elements of the
filtration  are defined by the images of $R^*(\overline{M}_{g,n})$
in the associated quotient sequence:
\begin{equation}
\label{ttpp}
R^*(\overline{M}_{g,n}) \rarr R^*(M^c_{g,n}) \rarr R^*(M^{rt}_{g,n})
\rarr R^*(C_{g,n}) \rarr 0
\end{equation}
Remarkably, the tautological rings of the strata
are conjectured to resemble
cohomology rings of compact manifolds.

A finite dimension graded algebra $R$ is Gorenstein with socle in
degree $s$ if there exists an evaluation isomorphism,
$$\phi: R^s \stackrel{\sim}{\rarr} \Q,$$
for which the bilinear pairings induced by the ring product,
$$R^r \times R^{s-r} \rarr R^s \stackrel{\phi}{\rarr} \Q,$$
are nondegenerate. The cohomology rings
of compact manifolds are
Gorenstein algebras.

\begin{conj}
The tautological rings of the filtration of $\overline{M}_{g,n}$
are finite dimensional Gorenstein algebras.
\end{conj}

The Gorenstein structure of $R^*(M_g)$ with socle in degree $g-2$
was discovered by Faber in his study
of the Chow rings of ${M_g}$ in low genus. The general conjecture
is primarily motivated by Faber's original work and can be found in
various stages in \cite{Fa}, \cite{hl}, and \cite{FP3}.

The application of the conjecture to the stratum
$C_{g,n}$ takes a special form
due to the choice of the underlying curve $C_g$. The conjecture
is stated for a nonsingular
curve $C_g$ defined over $\overline{\Q}$ or, alternatively,
for the tautological ring in $H^*(C_{g,n}, \Q)$.
The tautological ring of $C_{g,n}$ in Chow  is {\em  not} Gorenstein
for all $C_g$ by recent results of Green and Griffiths.

Two main questions immediately arise if the tautological rings are
Gorenstein algebras:
\begin{enumerate}
\item[(i)]  Can the ring structure be described explicitly?
\item[(ii)] Are the tautological rings associated to embedded compact manifolds
in the moduli space of pointed curves?
\end{enumerate}
 The tautological ring structures are implicitly determined
by the conjectural Gorenstein property and the Virasoro constraints \cite{GeP}.

As the moduli space of curves may be viewed as a special case of
the moduli space of maps,  a development of
these ideas may perhaps be pursued more
fully in Gromov-Witten theory. It is possible to
define a tautological ring for $\overline{M}_{g,n}(X,\beta)$ in the
context of the virtual class by {\em assuming} the Gorenstein property,
but no structure has been yet been conjectured.
Again, the Virasoro constraints in principle determine the tautological
rings.

\section{BPS states for threefolds} \label{bps}
\setzero\vskip-5mm \hspace{5mm }

Let $X$ be a  nonsingular projective variety over $\com$ of
dimension 3. Let $\{ \gamma_a \}_{a\in A} $ be a basis of $H^*(X,
\Z)$ modulo torsion. Let $\{ \gamma_a \}_{a\in D_2}$ and
$\{\gamma_a \}_{a\in D_{>2}} $ denote the classes of degree 2 and
degree greater than 2 respectively. The Gromov-Witten invariants
of $X$ are defined by integration over the moduli space of stable
maps (against the virtual fundamental class):
\begin{equation}
\label{klk}
\langle \gamma_{a_1}, \ldots, \gamma_{a_n} \rangle_{g,\beta}^X
 = \int_{[\overline{M}_{g,n}(X,\beta)]^{vir}}
\text{ev}_1^*(\gamma_{a_1}) \ldots \text{ev}_n^*(\gamma_{a_n}),
\end{equation}
where $\text{ev}_i$ is the $i$th evaluation map.
As the moduli spaces are Deligne-Mumford stacks, the
Gromov-Witten invariants take values in $\Q$.

Let $\{t_a \}$ be a set of variables corresponding to the classes $\{\gamma_a\}$.
The Gromov-Witten potential $F^X(t,\lambda)$ of $X$ may be written,
\begin{equation}
\label{gentt}
F^X = F^X_{\beta=0} +
\tilde{F}^X,
\end{equation}
as a sum of constant and nonconstant map contributions.

The
constant map contribution $F^X_{\beta=0}$ may be further divided by genus:
$$ F^X_{\beta=0} =
F^0_{\beta=0} + F^1_{\beta=0} + \sum_{g\geq 2} F^g_{\beta =0}.$$
The genus 0 constant contribution records the
classical intersection theory of $X$:
$$F^0_{\beta=0}= \lambda^{-2}
\sum_{a_1,a_2,a_3 \in A}  \frac{t_{a_3} t_{a_2} t_{a_1}}{3!}
\int_{X} \gamma_{a_1} \cup \gamma_{a_2} \cup \gamma_{a_3}.$$
The genus 1 constant contribution is obtained from a virtual
class calculation:
$$F^1_{\beta=0} = \sum_{a\in D_2} t_a \langle \gamma_a\rangle^X_{g=1, \beta=0} =
- \sum_{a\in D_2} \frac{t_a}{24}
\int_X \gamma_a \cup c_2(X).$$
Similarly, the genus $g\geq 2$
contributions are
$$F^g_{\beta=0} = \langle 1 \rangle^X_{g,\beta=0} =(-1)^g\frac{\lambda^{2g-2}}
{2}
 \int_X \left( c_3(X) -  c_1(X) \cup c_2(X)\right)          \cdot
\int_{\overline{M}_g} \lambda_{g-1}^3.$$
The Hodge integrals which arise here have been computed in \cite{FP}:
$$
\int_{\overline{M}_g}\lambda_{g-1}^3=\frac{|B_{2g}|}{2g}
\frac{|B_{2g-2}|}{2g-2}\frac1{(2g-2)!},
$$
where $B_{2g}$ and $B_{2g-2}$ are Bernoulli numbers.
The constant map contributions to $F^X$ are therefore completely
understood.

The second term in (\ref{gentt}) is  the
nonconstant map contribution:
$$\tilde{F}^X=
\sum_{g\geq 0} \sum_{\beta\neq 0} F^g_\beta.$$
Since the virtual dimension of the moduli space
$\overline{M}_g(X,\beta)$ is
$$\int_\beta c_1(X) + 3g-3+ 3-3g = \int_\beta c_1(X),$$
the classes $\beta$ satisfying $\int_{\beta} c_1(X)<0$
do not contribute to the potential $F^X$.
Therefore, $\tilde{F}^X$ may be divided into
two
sums:
\begin{eqnarray*}
\tilde{F}^X & = & \ \ \ \sum_{g\geq 0} \ \ \  \sum_{\beta\neq 0, \
\int_{\beta}c_1(X)=0} F^g_{\beta}
\\
& & + \sum_{g\geq 0} \ \ \ \sum_{\beta\neq 0,
\
\int_{\beta}c_1(X)>0} F^g_{\beta}.
\end{eqnarray*}

In case $\beta\neq 0$, we will write the series $F^g_{\beta}(t,\lambda)$
in the following
form:
$$F^g_{\beta} =\lambda^{2g-2}
 q^\beta
\sum_{n\geq 0} \frac{1}{n!} \sum_{a_1,\dots, a_n \in D_{>2}}
{t_{a_n} \cdots t_{a_1}}  \langle \gamma_{a_1} \cdots \gamma_{a_n}
\rangle^X_{g,\beta}.$$
The degree 2 variables $\{t_a\}_{a\in D_2}$ are
formally suppressed in $q$ via the divisor equation:
$$q^\beta= \prod_{a\in D_2} q_{a}^{\int_\beta \gamma_a}, \ \ \
q_a = e^{t_a}.$$
Cohomology classes of degree 0 or 1 do not appear in nonvanishing Gromov-Witten
invariants (\ref{klk}) for curve classes $\beta\neq 0$.

We will define new invariants $n_{\beta}^g(\gamma_{a_1}, \dots, \gamma_{a_n})$
for every genus $g$, curve class $\beta\neq 0$, and
classes $\gamma_{a_1}, \ldots, \gamma_{a_n}$.
Our primary interest will be in the case where the following
conditions hold:
\begin{enumerate}
\item[(i)] $\deg(\gamma_{a_i}) > 2$ for all $i$.
\item[(ii)]
$n+\int_{\beta}c_1(X) =\sum_{i=1}^n \deg(\gamma_{a_i}).$
\end{enumerate}
The invariants will be defined to satisfy the divisor equation
(which allows for the extraction of degree 2 classes $\gamma_a$)
and defined to vanish if degree 0 or 1 classes are inserted or if
condition (ii) is violated. If $\int_\beta c_1(X)=0$, then
$n_\beta^g$ is well-defined without cohomology insertions.

The new invariants $n_{\beta}^g(\gamma_{a_1}, \dots, \gamma_{a_n})$
are defined via Gromov-Witten theory by the following equation:
\begin{eqnarray*}
\tilde{F}^X & = &  \ \ \ \sum_{g\geq 0} \ \ \sum_{\beta\neq 0, \
\int_{\beta} c_1(X) =0} n_\beta^g \ \lambda^{2g-2} \sum_{d>0}
\frac{1}{d}\left( \frac{\sin
({d\lambda/2})}{\lambda/2}\right)^{2g-2} q^{d\beta} \\
& & + \sum_{g\geq 0} \ \ \sum_{\beta\neq 0, \ \int_{\beta} c_1(X)
>0} \ \ \sum_{n\geq 0} \frac{1}{n!}
 \sum_{a_1,\dots, a_n \in D_{>2}}{t_{a_n} \cdots t_{a_1}} \\
& & \ \ \ \ \ \ \ \ \ \ \cdot  \ n^g_\beta(\gamma_{a_1}, \cdots,
\gamma_{a_n})\
 \lambda^{2g-2} \left( \frac{\sin
({\lambda/2})}{\lambda/2}\right)^{2g-2+\int_{\beta}c_1(X)}
q^{\beta}.
\end{eqnarray*}
The above equation uniquely determines the invariants
$n^g_\beta(\gamma_{a_1}, \ldots, \gamma_{a_n})$.

\begin{conj}
For all nonsingular projective threefolds $X$,
\begin{enumerate}
\item[(i)]
the invariants $n_\beta^g(\gamma_{a_1}, \ldots, \gamma_{a_n})$ are
integers,
\item[(ii)] for fixed $\beta$, the invariants $n_\beta^g(\gamma_{a_1}, \ldots,
\gamma_{a_n})$ vanish for all sufficiently large genera $g$.
\end{enumerate}
\end{conj}

If $X$ is a Calabi-Yau threefold, the Gopakumar-Vafa conjecture is
recovered \cite{GV1}, \cite{GV2}.
Here,
the invariants $n_\beta^g$ arise as BPS state
counts in a study of Type IIA string theory on $X$
via M-theory. The outcome is a physical deduction of the
conjecture in the Calabi-Yau case.

Gopakumar and Vafa further propose a mathematical construction of
the Calabi-Yau invariants $n_\beta^g$ using moduli spaces of
sheaves on $X$. The invariants $n_\beta^g$ should arise as
multiplicities of special representations of ${\mathfrak{sl}}_2$
in the cohomology of the moduli space of sheaves. The local
Calabi-Yau threefold consisting of a curve $C$ together with a
rank 2 normal bundle $N$ satisfying $c_1(N)=\omega_C$ should be
the most basic case. Here the BPS states $n^g_d$ should be found
in the cohomology of an appropriate moduli space of rank $d$
bundles on $C$. A mathematical development of the proposed
connection between integrals over the moduli of stable maps and
the cohomology of the moduli of sheaves has not been completed.
However, evidence for the program can be found both in local and
global calculations in several cases \cite{bry}, \cite{hst},
\cite{kkv}.

The conjecture for arbitrary threefolds is motivated by
the Calabi-Yau case together with
the degeneracy calculations of \cite{P2}. Evidence
can be found, for example, in the low genus enumerative geometry
of $\proj^3$ \cite{Get}, \cite{P2}. If the conjecture is true, the invariants $n^g_\beta(\gamma_{a_1},
\ldots, \gamma_{a_n})$ of $\proj^3$ may be viewed as {\em defining} an integral enumerative
geometry of space curves for all $g$ and $\beta$. Classically the
enumerative geometry of space curves does not admit a uniform
description.

The conjecture does not determine the Gromov-Witten invariants
of threefolds. A basic related question is to find some means to
calculate higher genus invariants of Calabi-Yau threefolds. The basic
test case is
the quintic hypersurface in $\proj^4$. There are several approaches
to the genus 0 invariants of the quintic: Mirror symmetry, localization,
degeneration, and Grothendieck-Riemann-Roch \cite{can},\cite{giv2}, \cite{gath},\cite{giv1},\cite{yau}.
 But,
the higher genus invariants of the quintic are still beyond
current string theoretic and geometric techniques. The best tool
for the higher genus Calabi-Yau case, the holomorphic anomaly
equation,
 is not well understood in mathematics.
On the other hand, all the invariants of $\proj^3$ may be
in principle calculated by virtual localization \cite{GrP}.

\section{Virasoro constraints} \label{vir}
\setzero\vskip-5mm \hspace{5mm }

Let $X$ be a nonsingular projective variety over $\com$ of
dimension $r$. Let $\{ \gamma_a \}$ be a basis of $H^*(X,\C)$
homogeneous with respect to the Hodge decomposition,
$$\gamma_a \in H^{p_a,q_a}(X, \com),$$
The descendent
Gromov-Witten invariants of $X$ are:
$$\langle \tau_{k_1}(\gamma_{a_1}) \ldots \tau_{k_n}(\gamma_{a_n})
\rangle_{g,\beta}^X = \int_{[\overline{M}_{g,n}(X, \beta)]^{vir}}
\psi_1^{k_1} \ev_1^*(\gamma_{a_1}) \ldots  \psi_n^{k_n}
\ev_n^*(\gamma_{a_n}).$$ Let $\{t^a_k \}$ be a set of variables.
Let $F^X(t,\lambda)$  be the generating function of the descendent
invariants:
\begin{equation*}
F^X = \sum_{g\geq 0} \lambda^{2g-2} \sum_{\beta \in H_2(X,\Z)}
q^\beta \sum_{n\ge 0} \frac{1}{n!} \sum_{\substack{a_1\dots a_n \\
k_1 \dots k_n}} t_{k_n}^{a_n} \dots t_{k_1}^{a_1} \langle
\tau_{k_1}(\gamma_{a_1}) \dots \tau_{k_n}(\gamma_{a_n})
\rangle_{g,\beta}^X .
\end{equation*}
The partition function $Z^X$ is formed by exponentiating $F^X$:
\begin{equation} \label{GW}
Z^X = \exp(F^X).
\end{equation}

We will now define formal differential operators $\{ L_k\}_{k\geq
-1}$ in the variables $t^a_k$ satisfying the Virasoro bracket,
$$
[L_k,L_\ell] = (k-\ell) L_{k+\ell}.
$$
The definitions of the operators $L_k$ will depend only upon
the
following three structures of $H^*(X,\com)$:
\begin{enumerate}
\item[(i)] the intersection pairing $g_{ab} = \int_X \gamma_a \cup \gamma_b$,
\item[(ii)] the Hodge decomposition $\gamma_a \in H^{p_a,q_a}(X,\com)$,
\item[(iii)] the action of the anticanonical class $c_1(X)$.
\end{enumerate}
The formulas for the operators $L_k$ are:
\begin{align*}
L_k& = & \sum_{m=0}^\infty \sum_{i=0}^{k+1} &  \Bigl( [b_a\!+\!m]^k_i (C^i)^b_a
\tilde{t}^a_m \p_{b,m+k-i} \\
& & &  + \frac{\hbar}{2} (-1)^{m+1}
[b_a\!-\!m\!-\!1]^k_i (C^i)^{ab} \p_{a,m} \p_{b,k-m-i-1} \Bigr) \\ & & &
{\hspace{-40pt}}
+ \frac{\lambda^{-2}}{2} (C^{k+1})_{ab} t^a_0 t^b_0 \\ & & &
{\hspace{-40pt}}
+ \frac{\delta_{k0}}{48}
\int_X \bigl( (3-r)c_r(X)-2c_1(X)c_{r-1}(X) \bigr) ,
\end{align*}
where the Einstein convention for summing over the repeated indices
$a,b \in A$ is followed.

Several terms require definitions.
For each class $\gamma_a$, a half integer $b_a$ is obtained
from the Hodge decomposition,
$$b_a  =p_a+(1-r)/2.$$
The combinatorial factor $[x]^k_i$ is defined by:
$$
[x]^k_i = e_{k+1-i}(x,x+1,\dots,x+k) ,
$$
where $e_k$ is the $k$th elementary symmetric function.
The matrix $C_a^b$ is determined by the action of the anticanonical
class,
$$C_a^b \gamma_b = c_1(X) \cup \gamma_a.$$
The indices of $C$ are lowered and raised by the metric $g_{ab}$ and its
inverse $g^{ab}$.
The terms $\tilde{t}^a_m$ and $\p_{a,m}$ are  defined by:
\begin{eqnarray*}
\tilde{t}^a_m & = & t^a_m-\delta_{a0}\delta_{m1},  \\
 \p_{a,m}&  = & \p/\p t^a_m,
\end{eqnarray*}
where both are understood to vanish if
$m<0$.

\begin{conj}
For all nonsingular projective varieties $X$, $L_k(Z^X)=0$.
\end{conj}

The conjecture for varieties $X$ with only $(p,p)$ cohomology was
made by Eguchi, Hori, and Xiong \cite{EgHX}. The full conjecture
involves ideas of Katz. In case $X$ is a point, the constraints
specialize to the known Virasoro formulation of Witten's
conjecture \cite{K}, \cite{W} (see also \cite{big}). After the
point, the simplest varieties occur in  two basic families: curves
$C_g$  of genus $g$ and projective spaces $\proj^n$ of dimension
$n$. A proof of the Virasoro constraints for target curves $C_g$
is presented in a sequence of papers \cite{OP1}, \cite{OP2},
\cite{OP3}. Givental has recently proven the Virasoro constraints
for the projective spaces $\proj^n$ \cite{G1}, \cite{G2},
\cite{G3}. The two families of varieties are quite different in
flavor. Curves are of dimension 1, but have non-$(p,p)$
cohomology, non-semisimple quantum cohomology, and do not always
carry torus actions. Projective spaces cover all target
dimensions, but have algebraic cohomology, semisimple quantum
cohomology, and always carry torus actions.

The Virasoro constraints are especially appealing from the
point of view of algebraic geometry as
 {\em all} nonsingular projective varieties are covered.
While many aspects of Gromov-Witten theory may be more naturally
pursued in the symplectic category,
the Virasoro constraints appear to require more than a symplectic
structure to define. For example, the bracket
$$[L_1, L_{-1}] = 2 L_0,$$
depends upon formulas expressing the Chern numbers,
$$\int_X c_r(X), \ \ \int_X c_1(X)c_{r-1}(X),$$
in terms of the Hodge numbers $h^{p,q}$ of $X$ (see \cite{lib}]).

The Virasoro constraints may be a shadow of a deeper connection
between the Gromov-Witten theory of algebraic varieties and
integrable systems. In case the target is the point or the
projective line, precise connections have been made to the KdV and
Toda hierarchies respectively. The connections are proven by
explicit formulas for the descendent invariants in terms of matrix
integrals (for the point) and vacuum expectation in $\LV$ (for the
projective line) \cite{K}, \cite{big}, \cite{OP2}. The extent of
the relationship between Gromov-Witten theory and integrable
systems is not known. In particular, an understanding of the
surface case would be of great interest. Perhaps a link to
integrable systems can be found in the circle of ideas involving
Hilbert schemes of points, Heisenberg algebras, and G\"ottsche's
conjectures concerning the enumerative geometry of linear series.

Finally, one might expect Virasoro constraints to hold in the
context of Gromov-Witten theory relative to divisors in the target
$X$. For the relative theory of 1-dimensional targets $X$,
Virasoro constraints have been found and play a crucial role in
the proof of the Virasoro constraints for the absolute theory of
$X$ \cite{OP3}.

\end{document}